\documentclass[11pt]{amsart}
\usepackage{amsfonts,amssymb}
\usepackage{amscd}
\usepackage[colorlinks=true]{hyperref}
\usepackage{hyperref}
\hypersetup{colorlinks,citecolor=blue}

\newcommand{\BZ}{{\mathbb{Z}}}
\newcommand{\BN}{{\mathbb{N}}}
\newcommand{\BR}{{\mathbb{R}}}

\newcommand{\gD}{\Delta}
\newcommand{\gd}{\delta}
\newcommand{\gb}{\beta}

\newcommand{\gC}{\Gamma}
\newcommand{\gc}{\gamma}
\newcommand{\gs}{\sigma}
\newcommand{\gS}{\Sigma}
\newcommand{\gO}{\Omega}
\newcommand{\go}{\omega}
\newcommand{\gep}{\epsilon}
\newcommand{\gl}{\lambda}
\newcommand{\ga}{\alpha}
\newcommand{\gt}{\tau}

\newcommand{\gph}{\varphi}

\newcommand{\ti}[1]{\tilde{#1}}

\newtheorem{prop}{Proposition}[section]
\newtheorem{thm}[prop]{Theorem}
\newtheorem{lem}[prop]{Lemma}
\newtheorem{cor}[prop]{Corollary}

\theoremstyle{definition}
\newtheorem{Ack}[prop]{Acknowledgments}

\newtheorem{defn}[prop]{Definition}
\newtheorem{rem}[prop]{Remark}

\begin{document}
\author{T. Gelander, A. Karlsson, G.A. Margulis}

\thanks{T.G. partially supported by NSF grant
DMS-0404557. A.K. partially supported by VR grant 2002-4771. G.A.M
partially supported by NSF grant DMS-0244406. T.G. and G.A.M partially supported by BSF grant 2004010. }

\date{\today}

\title{Superrigidity, generalized harmonic maps and uniformly convex spaces}
\maketitle


\medskip

\begin{abstract}
We prove several superrigidity results for isometric actions on
Busemann non-positively curved uniformly convex metric spaces. In particular we generalize
some recent theorems of N. Monod on uniform and certain
non-uniform irreducible lattices in products of locally compact
groups, and we give a proof of an unpublished result on
commensurability superrigidity due to G.A. Margulis. The proofs rely on
certain notions of harmonic maps and the study of their existence,
uniqueness, and continuity.
\end{abstract}

Ever since the first superrigidity theorem for linear
representations of irreducible lattices in higher rank semisimple
Lie groups was proved by Margulis in the early 1970s, see
\cite{Margulis-sr} or \cite{Margulisb}, many extensions and
generalizations were established by various authors, see for
example the exposition and bibliography of \cite{Jost} as well as
\cite{Pansu}.
A superrigidity statement can be read as follows:\\
Let
\begin{itemize}
\item $G$ be a locally compact group,

\item $\gC$ a subgroup of $G$,

\item $H$ another locally compact group, and

\item $f:\gC\to H$ a homomorphism.
\end{itemize}
Then, under some certain conditions on $G,\gC,H$ and $f,$ the
homomorphism $f$ extends uniquely to a continuous homomorphism
$F:G\to H$. In case $H=\text{Isom}(X)$ is the group of isometries
of some metric space $X$, the conditions on $H$ and $f$ can be
formulated in terms of $X$ and the action of $\gC$ on $X$.

In the original superrigidity theorem \cite{Margulis-1} it was assumed that $G$ is a
semisimple Lie group of real rank at least
two\footnote{Superrigidity theorems were proved later also for
lattices in the rank one Lie groups $\text{SP}(n,1), F_4^{-20}$
see \cite{Corlette} and \cite{GromovSchoen}. It seems however that
the same phenomenon holds in these cases for different reasons.}
and $\gC\leq G$ is an irreducible lattice. It is not clear how to
define a rank for a general topological group. One natural
extension, although not a generalization, of the notion of higher
rank is the assumption that $G$ is a non-trivial product.
Margulis \cite{Margulis-1} also proved a superrigidity theorem for commensurability subgroups in semisimple Lie groups. 
The target in these superrigidity theorems was the 
group of isometries of a Riemannian symmetric space of non-compact type or an affine building.

It was later realized in an unpublished manuscript of
Margulis \cite{Margulis} which was circulated in the 1990s (cf.
\cite{Jost}), that superrigidity for commensurability subgroups
extends to a very general setting: a general locally compact,
compactly generated $G$
and a target group being the isometry group of a complete
Busemann non-positively curved uniformly convex metric space.

In this paper we establish quite general superrigidity theorems
for actions of irreducible lattices in products of locally compact
groups on Busemann non-positively curved uniformly convex metric spaces.
We also include a proof of the unpublished result for
commensurability subgroups mentioned above, since the methods are
similar.

Our method relies on certain notions of generalized harmonic
maps. The main part is the proof of their existence (Theorem \ref{thm:existence}, Proposition \ref{prop:existcont}) which is of independent interest and may have other applications. Once the existence is established the superrigidity results follow from the nice properties of these maps.

Our results for lattices in products generalize recent theorems of N. Monod
\cite{Monod} for actions on CAT(0) spaces. Loosely speaking, the argument of Monod \cite{Monod} is divided into three steps: 1. Inducing the lattice action to an action of the ambient group on the space of square integrable equivariant maps. 2. Proving a splitting theorem for actions of product groups on CAT(0) spaces (generalizing earlier results for Riemannian manifolds and general proper CAT(0) spaces, c.f. \cite{Br-Ha} p. 239). 3. Using the splitting of the induced space to obtain an invariant subset of the original space on which the lattice action extends.
The proof we give here is in a sense more direct and therefore applies in a more general setup where the splitting result does not hold.
The spaces considered in the current paper, namely Busemann non-positively curved uniformly convex metric spaces, generalize CAT(0) spaces in a similar manner as (uniformly convex) Banach spaces generalize Hilbert spaces and (uniformly convex NPC) Finsler manifolds generalize NPC Riemannian manifolds. In particular, CAT(0) spaces are BNPC and UC, but these conditions are much weaker than CAT(0), for instance Hilbert spaces, while being the least convex among CAT(0) spaces, are the most convex among UC Banach spaces (We refer the reader to \cite{harpe},\cite{Ne},\cite{Up} for examples).
There are several technical difficulties that could be avoided by requiring stronger assumptions on the spaces considered (for example, Proposition \ref{prop:Lp} with $p=2$ is obvious for CAT(0) spaces). In particular, when assuming the CAT(0) condition, our proof is significantly simplified.
Note that we also obtain some new results for actions on CAT(0) spaces (c.f. Theorem \ref{thm3}).

Sections 1 through 7 deal with lattices in products, and
Section 8 deals with commensurability subgroups. The argument in Section 8 is slightly simpler and more detailed than the original proof given in \cite{Margulis}, however we assume here a slightly stronger convexity assumption on the space $X$.

\begin{Ack} We would like to thank Nicolas Monod for several
helpful discussions, and the anonymous referees for their remarks and suggestions.
\end{Ack}


\section{Assumptions and conclusion}\label{sec:statements}
We shall first discuss some properties of the groups, the spaces
and the actions under consideration, and then state our main
results.

By saying that a lattice $\gC$ in a locally compact group with a
specified decomposition $G=G_1\times\ldots G_n$ is {\it
irreducible} we mean that $\gC (\prod_{j\neq i}G_j)$ is
dense\footnote{Note that the projection of $\gC$ to a sub-product
$\prod_{j\in J}G_j$ with $|J|>1$ is not necessarily dense.} in $G$
for each $1\leq i\leq n$.


A complete geodesic metric space $X$ is said to be {\it Busemann
non-positively curved}, or shortly BNPC, if the distance between
any two constant speed geodesics is a convex function. In
particular, $X$ is a uniquely geodesic space, \noindent i.e. any
two points $x,y\in X$ are joined by a unique arc whose length is
$d(x,y)$. We shall denote the midpoint of $x$ and $y$ by
$\frac{x+y}{2}$.

A uniquely geodesic metric space $X$ is said to be {\it strictly convex} if
$d(x,\frac{y_1+y_2}{2})<\max\{ d(x,y_1),d(x,y_2)\}$,
$\forall x,y_1,y_2\in X$ with $y_1\neq y_2$. We shall say that $X$ is {\it weakly uniformly convex} (or shortly WUC) if
additionally for any $x\in X$,
the modulus of convexity function
$$
 \gd_{x}(\gep,r):=\inf\{r-d(x,\frac{y_1+y_2}{2}):y_i\in X,d(x,y_i)\leq r,d(y_1,y_2)\geq\gep
 r\}
$$
is positive for any $\gep,r>0$. We shall say that $X$ is {\it uniformly convex} (or shortly UC) if 
$$
 \forall\gep> 0, \exists \gd (\gep)>0~\text{such that}~\forall r>0,x\in X,~\gd_x(\gep,r)\geq\gd(\gep)\cdot r.
$$
We do not know an example of a BNPC WUC space which is not UC, and it is
conceivable that BNPC and WUC imply UC. We will allow ourselves to assume whenever we find it convenient that the metric space under consideration is UC, although some of our results can be proved under weaker convexity assumptions. For instance Theorem \ref{thmcommensurator} was proved in \cite{Margulis} under the assumption that $X$ is BNPC and WUC. All the main results stated below remains true if one assume only that $X$ is uniformly convex with respect to some point $x_0\in X$, i.e. that $\inf_{r>0}\frac{1}{r}\gd_{x_0}(\gep,r)>0,~\forall\gep>0$ (and the proofs require only minor changes).

The \emph{projection} of a point $x\in X$ to a closed convex set $C\subset X$ is
the point $p\in C$ closest to $x$. It is a standrd fact that
projections exists uniquely in spaces which are weakly uniformly
convex and complete.

Suppose that $\gC$ acts on $X$ by isometries. The action is called
{\it $C$--minimal} if there is no non-empty closed convex proper
$\gC$--invariant subset of $X$.

For any subset $\gS\subset\gC$ we associate the displacement
function:
$$
 d_{\gS}(x):=\max_{\gs\in\gS}d(\gs\cdot x,x).
$$
We shall say that {\it $d_\gS$ goes to infinity} and shortly write
$d_\gS\to\infty$ if $\lim_{x\to\infty}d_\gS (x)=\infty$, where
$x\to\infty$ means that $x$ eventually gets out of any ball.

Following \cite{Monod}, we shall say that an action
$\gC\circlearrowleft X$ is {\it reduced} if there is no unbounded
closed convex proper subset $Y\subset X$ which is of finite
Hausdorff distance from $\gc\cdot Y$ for any $\gc\in\gC$.

Recall that a {\it Clifford isometry} $T$ of a metric space $X$ is
a surjective isometry $T:X\to X$ for which $d(x,T(x))$ is constant
on $X$. A necessary assumption in some of the superrigidity theorems below is that $X$ has no
non-trivial Clifford isometries. However, for typical spaces this
assumption follows from reduceness of $\text{Isom}(X)$ (c.f.
Corollary \ref{R->NCL}).

The assumption that the action is reduced is very strong, and by
requiring it we can prove a quite general statement.

\begin{thm}\label{thm1}
Let $G=G_1\times\ldots\times G_n,~n>1$ be a locally compact
compactly generated topological group, and $\gC$ an irreducible
uniform lattice in $G$. Let $X$ be a complete Busemann
non-positively curved uniformly convex metric space without
non-trivial Clifford isometries. Assume that $\gC$ acts on $X$ by
isometries and that this action is reduced and has no global fixed
point. Then the action of $\gC$ extends uniquely to a continuous
$G$--action, and this $G$--action factors through one of the $G_i$'s.
\end{thm}

\begin{rem}\label{rem:weaker-assumption}
$(i)$ In Section \ref{sec:non-uniform} we shall give a
generalization of Theorem \ref{thm1} for non-uniform lattices
which are $p$--integrable and weakly cocompact.

$(ii)$ The theorem also remains true when $G$ is not assumed to be
compactly generated but only a $\gs$--compact group\footnote{One
should only choose the function $h$ in the definition of the
energy $E$ below more carefully so that the energy of some
$\gC$--equivariant map (in $L_2(\gC\backslash G,X)$) would be
finite.} if we add the assumption that $d_\gS\to\infty$ for some
finite subset $\gS\subset\gC$. Note however that when $G$ is not
compactly generated, $\gC$ is not finitely generated. Hence this
additional assumption is perhaps not very natural.
\end{rem}


%

Even when an action is reduced, it is not clear how to verify
this. Moreover, in many cases, the action can be extended under
weaker assumptions. By requiring stronger assumptions on the space
$X$ we can drop this assumption. The following two theorems imply
(and generalize) Margulis' original theorem in the case where $G$
is not simple. For an explanation why Theorem \ref{thm2} implies
Margulis Theorem we refer to \cite{Monod} where the same result is
proved under the assumption that $X$ is CAT(0). Note that Remark
\ref{rem:weaker-assumption} applies also for Theorems \ref{thm2}
and \ref{thm3}.

\begin{thm}\label{thm2}
Let $G,\gC$ be as in Theorem \ref{thm1}, let $X$ be a complete
BNPC uniformly convex space, and assume further that $X$ is proper
(i.e. the closed balls in $X$ are compact). Assume that $\gC$ acts
on $X$ without a global fixed point. Then there is a non empty
closed invariant subset $\mathcal{L}$ of the visual boundary $\partial X$ (c.f. \cite{Papado} for definition) on which
the $\gC$ action extends to a continuous $G$ action which factors
through some $G_i$.
\end{thm}

We will say that the space $X$ is {\it geodesically complete}
(resp. {\it uniquely geodesically complete}) if any geodesic
segment in $X$ is contained in a (unique) two sided infinite
geodesic.

\begin{thm}\label{thm3}
Let $G$ and $\gC$ be as above, and let $\gS$ be a finite
generating set of $\gC$. Let $X$ be a complete uniquely
geodesically complete CAT(0) space without Euclidean
factors\footnote{Recall that a CAT$(0)$ space $X$ admits a
non-trivial Clifford isometry iff it has an Euclidean factor, i.e.
can be decomposed as a direct product $X'\times\BR$ (c.f.
\cite{Br-Ha} p. 235).} with the additional assumption that if two
geodesic segments are parallel then the corresponding geodesics are parallel.
Suppose that $\gC$ acts on $X$ by isometries with $d_\gS\to\infty$. Then
there is a $\gC$--invariant geodesically complete closed subset
$Y\subset X$, such that the $\gC$--action on $Y$ extends to a
continuous $G$--action. Moreover, $Y$ can be decomposed to an
invariant direct product $Y=Y_1\times\ldots\times Y_n$ such that
each $G_i$ acts trivially on each $Y_j$ with $j\neq i$.
\end{thm}

\begin{rem}
(1) If $X$ is CAT(0) and the distance function between any two
geodesics is analytic (except at intersection points) then all
assumptions on $X$  (in \ref{thm3}) are satisfied.
(2) Trees and buildings are usually not uniquely geodesically complete.
\end{rem}


\section{Some remarks about the assumptions}

\subsection{Parallel segments}
We shall say that two segments $[a,b],[x,y]\subset X$
are parallel, and write $[a,b]~\|~[x,y]$, if
$$
 d(a,x)=d(b,y)=d\big(\frac{a+b}{2},\frac{x+y}{2}\big).
$$

Recall the following fact:

\begin{lem}[Busemann \cite{Busemann}, Th. 3.14]\label{lem:parallelogram}
Whenever $X$ is BNPC
$$
 [a,b]~\|~[x,y]\Leftrightarrow [a,x]~\|~[b,y].
$$
\end{lem}

\subsection{Intersection property for convex sets}

A real valued function on $X$ is called {\it convex} if it is convex
in the usual sense when restricted to any geodesic segment, i.e. if its sub-level sets are convex. 

\begin{lem}\label{lem-bounded-convex}
Let $X$ be a weakly uniformly convex complete metric space. Then
any collection of closed bounded convex sets has the finite intersection property. If $f$ is a
convex function on $X$ which satisfies $f(x)\to\infty$ when
$x\to\infty$, then $f$ attains a minimum in $X$.
\end{lem}

\begin{proof}
To prove the first claim let $C_\ga$ be a descending net of closed
convex sets which are contained in some ball. Fix $x_0\in X$ and
let $x_\ga$ be the projection of $x_0$ to $C_\ga$. Then
$d(x_0,x_\ga)$ is a non-decreasing bounded net of non-negative
numbers and hence has a limit. If this limit is $0$ then $x_0$
belongs to the intersection $\cap C_\ga$, and if it is positive
one shows that $x_\ga$ must be a Cauchy net as follows: if
$\gb>\ga$ then $x_\gb$ belongs to $C_\ga$ and if $\ga$
is "large" then $d(x_0,x_\gb)$ is "almost" the same as $d(x_0,x_\ga)$ and hence $d(x_\gb,x_\ga)$ must be "small" for otherwise
$\frac{x_\ga+x_\gb}{2}$ would be closer to $x_0$ than $x_\ga$, by
positivity of the modulus of convexity function. Clearly $\lim
x_\ga\in\cap C_\ga$, and hence $\cap C_\ga\neq\emptyset$.

The second claim follows from the first by taking the convex sets
to be non-empty sub-level sets of $f$.
\end{proof}

\subsection{Linear growth of convex functions}
The following lemma will be used in the proof of the existence of
a harmonic map.

\begin{lem}\label{lem:b>0}
Let $X$ be a geodesic metric space and let $f:X\to\BR$ be a convex function.
Assume that $f(x)\to\infty$ when $x\to\infty$. Then there exists
$b>0$ such that $f(x)\geq b\cdot d(x,x_0)-\frac{1}{b},~\forall
x\in X$. Moreover if $X$ is weakly uniformly convex and
$f(x)>0,\forall x\in X$, then there exist $b>0$ such that
$f(x)>b\cdot d(x_0,x),\forall x\in X$.
\end{lem}

\begin{proof}
The second statement follows from the first using Lemma \ref{lem-bounded-convex}.
Now, assuming the contrary, there must be a sequence $x_n$ in $X$ such that
\begin{itemize}
\item $f(x_n)\leq\frac{1}{n}d(x_n,x_0)$, and

\item $d(x_n,x_0)\geq n^2$.
\end{itemize}
Take $y_n$ to be the point of distance $n$ from $x_0$ on the
geodesic segment $[x_0,x_n]$. Then $y_n\to\infty$ but
$$
 f(y_n)\leq
 \frac{d(y_n,x_n)f(x_0)+d(y_n,x_0)f(x_n)}{d(x_0,x_n)}\leq
 f(x_0)+\frac{\frac{d(y_n,x_0)}{n}d(x_0,x_n)}{d(x_0,x_n)}=f(x_0)+1.
$$
This however contradicts the assumption that
$f(x)\to\infty$ when $x\to\infty$.
\end{proof}

\subsection{Spaces with Clifford isometries}
For spaces with Clifford isometries, the analog of
Theorem \ref{thm1} is not true as shown by the simple example from
\cite{Monod} where $G$ is the discrete group $(\BZ/(2)\ltimes
\BZ)\times (\BZ/(2)\ltimes \BZ)$ and $\gC$ the index 2 irreducible subgroup
$\BZ/(2)\ltimes (\BZ\oplus\BZ )$, and $\gC$ acts on $\BR$, each
$\BZ$ by translation and the order 2 element by reflection trough
$0$; This action does not extend to $G$.

The typical example of a space with many Clifford isometries is a
Banach space. A weaker superrigidity result for isometric and, more
generally, for uniformly bounded affine actions of irreducible
lattices on uniformly convex Banach spaces was proved in
\cite{BFGM} (cf. Theorem D there).

On the other hand, spaces with non-trivial Clifford isometries
admit a canonical non-trivial invariant metric foliation and one can prove
superrigidity theorems for the induced action on the space of
leaves. Moreover, an action on such a space cannot be reduced,
unless it is a Banach space. 

\begin{prop}\label{clifford-folliation}
Let $X$ be a complete BNPC metric space, and assume that the set $\text{CL}(X)$ of all Clifford isometries of $X$ form a subgroup of $\text{Isom}(X)$. Then $\text{CL}(X)$ is normal and abelian. The orbits of $\text{CL}(X)$ are all isometric to some fixed uniformly
convex real Banach space. The quotient space $X/\text{CL}(X)$ is
BNPC. The induced action of $\text{Isom}(X)$ on $X/\text{CL}(X)$
is by isometries. Furthermore, if $X$ is uniformly convex then so
is $X/\text{CL}(X)$.
\end{prop}

\begin{proof}
It follows from Lemma
\ref{lem:parallelogram} that $S^{-1}T^{-1}ST(x)=x,~\forall x\in X$, and
hence, since $S,T\in\text{CL}(X)$ are arbitrary, $\text{CL}(X)$ is abelian.
Clearly $\text{CL}(X)\triangleleft\text{Isom}(X)$. To give the
structure of a Banach space to $\text{CL}(X)\cdot x$ note that it
follows from our assumptions and Lemma \ref{lem:parallelogram} that
$T^{\frac{1}{2}}(x):=\frac{x+T(x)}{2}$ is also a Clifford isometry
and hence we can define a multiplication by a dyadic number, and
by continuity we can define a multiplication by any real number.
It is also easy to verify that $X/\text{CL}(X)$ is BNPC, and also
uniformly convex if $X$ is. Any isometry $T$ induces a
1--Lipschitz map on $X/\text{CL}(X)$, which forces it, as $T^{-1}$
is also a 1--Lipschitz, to be an isometry.
\end{proof}

\begin{rem}
The assumption that $\text{CL}(X)$ is a group holds in many cases, e.g. it holds for CAT(0) spaces (cf. \cite{Br-Ha} p. 235), as well as spaces for which parallelity of segments is a transitive relation, i.e. 
$([a_1,b_1]~\|~[a_2,b_2])\wedge ([a_2,b_2]~\|~[a_3,b_3])\Rightarrow [a_1,b_1]~\|~[a_3,b_3]$.
\end{rem}

\begin{cor}\label{R->NCL}
Let $X$ be a BNPC complete metric space which is not isometric to
a Banach space, and assume that $\text{CL}(X)$ is a group. Suppose that $\text{CL}(X)$ is non-trivial. Then the action
of $\text{Isom}(X)$ on $X$ is not reduced.
\end{cor}

\begin{proof}
Choose $y\in X$ and let $Y=\text{CL}(X)\cdot y$. Then $Y$ is a
closed convex unbounded proper subset which is equidistant to its
image under any isometry.
\end{proof}

Still, for spaces with Clifford isometries one can prove the following:

\begin{thm}\label{thm:CL}
Let $\gC$ and $G$ be as in Theorem \ref{thm1}, and let $X$ be a complete BNPC  uniformly convex metric space such that $\text{CL}(X)$ is a group.
Assume that $\gC$ acts on $X$ with $d_\gS\to\infty$ where $\gS\subset\gC$ is a finite generating set, and that the induced action of $\gC$ on $X/\text{CL}(X)$ is reduced. Then the $\gC$ action on $X/\text{CL}(X)$ extends uniquely to a continuous $G$ action which factors through some $G_i$.
\end{thm}

We shall not elaborate on the proof of Theorem \ref{thm:CL}, which requires only a small modification in Step (2) of the proof of Theorem \ref{thm1}.

\subsection{Reduced actions, displacement functions and
$C$--minimality}

Suppose that a group $\gC$ is generated by a finite set $\gS$ and
acts by isometries on a space $X$.

\begin{lem}\label{lem:reduceness-->displacement} 
Assume that $X$ is a BNPC space which is not isomorphic to a Banach
space and that $\text{CL}(X)$ is trivial. If the action is reduced then the displacement function
$d_\gS\to\infty$.
\end{lem}

\begin{proof}
Since $\text{CL}(X)$ is trivial, if $\gc\in\gS$ acts non-trivially on $X$, $d_\gc(x)=d(x,\gc\cdot x)$ is non-constant. This implies that
the function $d_\gS$ is non-constant. To see this one can argue by
contradiction as follows: Suppose $d_\gS(x)\equiv m$, let $x\in X$
be a point where $J(x)=|\{\gc\in\gS:d_\gc (x)=m\}|$ is minimal,
let $\gc\in\gS$ be such that $d_\gc (x)=m$ and let $y\in X$ be
another point where $d_\gc (y)<m$. Then $J(\frac{x+y}{2})\leq
J(x)-1$ in contrary to the minimality of $J(x)$.

Now, since $X$ is BNPC the function $d_\gS$ is convex, being the
maximum of the convex functions $\{d_\gc:\gc\in\gS\}$. Since the
action is reduced, each proper sub-level set $\{x\in X:d_\gS
(x)\leq a\},~a>0$, is bounded, being convex and of bounded Hausdorff distance from its translation by any $\gc\in\gC=\langle\gS\rangle$. On the other hand, one can show that a bounded
non-constant convex function on a unbounded geodesic space must have an unbounded proper sub-level
set. We may of course assume that $X$ is unbounded for otherwise there is nothing to prove. It follows that $d_\gS$ is unbounded with bounded sub-level
sets, i.e. that $d_\gS\to\infty$.
\end{proof}

\begin{rem}
In view of Corollary \ref{R->NCL}, the assumption that $\text{CL}(X)$ is trivial in Lemma \ref{lem:reduceness-->displacement} can be replaced by the assumption that it is a group.
\end{rem} 

\begin{lem}\label{lem:C-min}
Assume that $X$ is complete, BNPC and uniformly convex. If
$d_\gS\to\infty$ then there exists a minimal closed convex
invariant subset.
\end{lem}

\begin{proof}
We may assume that there is no global fixed point. Let $C_\ga$ be
a descending chain of closed convex invariant sets. We claim that
the intersection $\cap C_\ga$ is non-empty. By Lemma
\ref{lem-bounded-convex} it is enough to show that they all
intersect some given ball. Fix $x_0\in X$ and let $x_\ga$ be the
projection of $x_0$ to $C_\ga$. Let $b$ be the constant from Lemma
\ref{lem:b>0} applied to the convex function $d_\gS$ and let
$a=d_\gS (x_0)$. Since $X$ is uniformly convex there is $R>0$ such
that $\gd_{x_0} (b,r)>a$ for all $r\geq R$. We claim that
$d(x_0,x_\ga)\leq R$. Suppose in contrary that $d(x_0,x_\ga)>R$,
and let $\gc\in\gS$ be such that $d(\gc\cdot x_\ga,x_\ga)\geq
b\cdot d(x_0,x_\ga)$. Since $d(x_0,\gc\cdot x_\ga)\leq
a+d(x_0,x_\ga)$ this implies that $d(x_0,\frac{x_0+\gc\cdot
x_\ga}{2})<d(x_0,x_\ga)$, in contrary to the definition of
$x_\ga$.
Thus, we can apply Zorn lemma and conclude that there is a minimal
invariant closed convex set.
\end{proof}

From the existence and uniqueness of circumcenters in weakly uniformly
convex, complete spaces, it follows that a bounded minimal invariant convex set must be a point. Indeed, the projection of the circumcenter\footnote{Note that if $X$ is not CAT$(0)$, the circumcenter of a closed bounded convex set may lie outside the set.} of a closed bounded invariant convex set to the set is fixed point. 

Therefore:

\begin{cor}
Assume that $X$ is complete, BNPC and uniformly convex. If the
action of $\gC$ on $X$ is reduced, and has no global fixed points,
then it is also $C$--minimal.
\end{cor}

By the exact same argument as in the proof of Lemma
\ref{lem:C-min} one can show:

\begin{lem}\label{S-min}
Assume that $X$ is a complete, uniquely geodesically complete,
BNPC and uniformly convex. If $d_\gS\to\infty$ then there exists a
minimal closed geodesically complete invariant subset.
\end{lem}

\subsection{Spaces of $p$--integrable maps}\label{par:Lp}

Let $(\gO,\mu)$ be a probability space and let $1<p<\infty$. We denote by $L_p(\gO,X)$ the space of measurable maps
$\gph :\gO\to X$ which satisfy $\int_\gO d(\gph (w),x_0)^p<\infty$ for $x_0\in X$\footnote{Obviously, the finiteness of this integral is independent of the choice of $x_0$.} with the distance
$$
 \rho (\gph,\psi)=\big(\int_\gO d(\gph (w),\psi (w))^p\big)^\frac{1}{p}.
$$
We shall make use of the following proposition whose technical proof 
might be skipped at first reading (related material may be found in \cite{Foertsch}).

\begin{prop}\label{prop:Lp}
Suppose that $X$ is complete BNPC and uniformly convex. 
Then for any $1<p<\infty$, the space $L_p(\gO,X)$ is also complete BNPC and 
uniformly convex.
\end{prop} 

\begin{proof}
The completeness of $L_p(\gO,X)$ follows from that of $X$ by a straightforward argument. Similarly the BNPC property follows from that of $X$ since three points $f,g,h\in L_p(\gO,X)$ lie on a common geodesic iff $f(w),g(w)$ and $h(w)$ lie on a common geodesic in $X$ for almost every $w\in\gO$.

Let us show that $L_p(\gO,X)$ is uniformly convex.
Let $\gd(\gep)>0$ be the associated constant in a linear lower bound for the modulus of convexity $\gd_{x_0}(\gep,r)$. We may assume that $\gd(\gep/4)$ is sufficiently small comparing to $\gep$ to satisfy Inequality (\ref{eq:2}) as well as the last calculation in the proof. 
Let $\gb_p$ be the modulus of convexity function of the Banach space $L_p(0,1)$. Set
$$
 \gt (\gep)=\gb_p(\gd^4({\gep}/{4})).
$$
We claim that if $\gd(\frac{\gep}{4})$ is chosen small enough then for any $\psi,\gph_1,\gph_2\in L_p(\gO,X)$ which satisfy $\rho (\gph_i,\psi)\leq r,~i=1,2$ and $\rho (\gph_1,\gph_2)\geq\gep r$ we have $\rho (\frac{\gph_1+\gph_2}{2},\psi)\leq r(1-\gt (\gep ))$. To see this, note first that if
\begin{equation}\label{eq:1}
 \big(\int_\gO |d(\gph_1(w),\psi(w))-d(\gph_2(w),\psi(w))|^p\big)^\frac{1}{p}\geq\gd^4({\gep}/{4})r
 \end{equation}
then the this claim follows from the uniform convexity of $L_p(0,1)$ since by BNPC
$$
 d(\frac{\gph_1(w)+\gph_2(w)}{2},\psi(w))\leq \frac{d(\gph_1(w),\psi(w))+d(\gph_2(w),\psi(w))}{2},~\forall w\in\gO.
$$
We shall therefore assume below that Inequality (\ref{eq:1}) does not hold. Set
$$
 \gO':=\{ w\in\gO : |d(\gph_1(w),\psi(w))-d(\gph_2(w),\psi(w))|\geq \gd^2({\gep}/{4})d(\gph_1(w),\psi(w))\},
$$
then the negation of Inequality (\ref{eq:1}) implies
\begin{equation}\label{eq:a}
\big(\int_{\gO'}d(\gph_1(w),\psi(w))^p\big)^\frac{1}{p} <\gd^2({\gep}/{4})r,
\end{equation}
and hence, using the negation of (\ref{eq:1}) again, $\big(\int_{\gO'}d(\gph_2(w),\psi(w))^p\big)^\frac{1}{p} <(\gd^2(\frac{\gep}{4})+\gd^4(\frac{\gep}{4}))r$.
Thus
\begin{equation}\label{eq:2}
 \int_{\gO\setminus\gO'} d(\gph_1(w),\gph_2(w))^p\geq (\gep^p-(2\gd^{2}({\gep}/{4})+\gd^{4}({\gep}/{4}))^p)r^p\geq (\frac{\gep}{2}r)^p.
 \end{equation}
Let now
$$
 \gO''=\{ w\in\gO\setminus\gO': d(\gph_1(w),\gph_2(w))\geq \frac{\gep}{4}d(\gph_1(w),\psi(w))\}.
$$
Then from Inequality (\ref{eq:2}) it follows that
\begin{equation}\label{eq:b}
\big(\int_{\gO''}d(\gph_1(w),\psi(w))^p\big)^\frac{1}{p}\geq \frac{\gep}{4}r.
\end{equation}
Thus, by uniform convexity of $X$
$$
 \big(\int_{\gO''}d(\frac{\gph_1+\gph_2}{2},\psi)^p\big)^\frac{1}{p}\leq \big(\int_{\gO''}d(\gph_1,\psi)^p\big)^\frac{1}{p}-
 (\frac{\gep}{4}\gd ({\gep}/{4})-\gd^4({\gep}/{4}))r.
$$
Therefore, assuming $\gd (\frac{\gep}{4})$ is small enough, we have
$$
 \int_\gO \Big(d(\gph_1,\psi)^p-d(\frac{\gph_1+\gph_2}{2},\psi)^p\Big)=\big(\int_{\gO''}+ 
 \int_{\gO\setminus\gO''} \big)\big(d(\gph_1,\psi)^p-d(\frac{\gph_1+\gph_2}{2},\psi)^p\big)
$$
$$
 \geq\int_{\gO''}\Big(d(\gph_1,\psi)^p-d(\frac{\gph_1+\gph_2}{2},\psi)^p\Big)-\int_{\gO\setminus\gO''}
 \big(d(\gph_2,\psi)^p-d(\gph_1,\psi)^p\big)\vee 0
$$
$$
 \geq \big(1-[1-\frac{\gep}{4}\gd(\gep/4)+\gd^4(\gep/4)]^p\big)(\frac{\gep}{4}r)^p-p\gd^4({\gep}/{4})r^p
 \geq \gd^2({\gep}/{4})r^p.
$$
It follows that $\big(\int_\gO d(\frac{\gph_1+\gph_2}{2},\psi)^p\big)^\frac{1}{p}\leq (1-\gd^{2}(\frac{\gep}{4}))^\frac{1}{p}r$. Finally, if $\gd(\frac{\gep}{4})$ is sufficiently small then $(1-\gd^{2}(\frac{\gep}{4}))^\frac{1}{p}\leq 1-\gd^{4}(\frac{\gep}{4})$, and since $\gb_p(\gd^4(\frac{\gep}{4}))\leq \gd^4(\frac{\gep}{4})$ this completes the proof.
\end{proof}

Arguing similarly, and avoiding the part of $\gO$ where both $\gph_i$ are close to $\psi$ one can show: 

\begin{prop}\label{prop:Lp-UCld} 
If $X$ is uniformly convex for large distances, i.e. for some (and hence any) $x\in X$, $\lim_{r\to\infty}\frac{1}{r}\gd (\gep,r)>0$ for any $\gep>0$, then so is $L_p(\gO,X)$.
\end{prop}

Since we will not make use of this variant of Proposition \ref{prop:Lp}, we shall not elaborate on its proofs.

\begin{rem}
It is straightforward that if $X$ is CAT(0) then so is $L_2(\gO,X)$. Thus for CAT(0) spaces one can avoid the technical Proposition \ref{prop:Lp}.
\end{rem}


\section{Definition and existence of generalized harmonic maps}\label{sec:existence}

Let $G=G_1\times G_2$ be a compactly generated locally compact
group and $\gC\leq G$ a uniform lattice. Let $\gO$ be a relatively
compact right fundamental domain for $\gC$, i.e.
$G=\bigsqcup_{\gc\in\gC}\gc\cdot\gO$. Assume that $X$ is a
complete BNPC uniformly convex metric space, and that $\gC$ acts
by isometries on $X$ with $d_\gS\to\infty$, where $\gS$ is a
finite generating set for $\gC$ containing the
identity\footnote{Note that in this section, we neither require
that the action $\gC\circlearrowleft X$ is reduced, nor that
$\text{CL}(X)$ is trivial, nor that $\gC\leq G$ is irreducible.}.

A function $\varphi:G\to X$ is said to be {\it $\gC$--equivariant}
if $\varphi (\gc g)=\gc\cdot\varphi (g)$ for any $\gc\in\gC,g\in
G$. Since such a function is determined by its restriction to
$\gO$ we will abuse notation and make no distinction between
$\gC$--equivariant functions and their restriction to $\gO$. In particular, by $L_2(\gO,X)$ we mean the space of $\gC$--equivariant measurable maps $\gO\to X$ whose restriction to $\gO$ is square integrable (see \ref{par:Lp}).

Fix a compact generating set $K$ of $G_1$ and define $h:G_1\to\BR$
by $h(g_1):=e^{-|g_1|_K^2+1}$ where $|~|_K$ is the word norm
associated to $K$, i.e. $|g_1|_K=\min\{ k:K^k\ni g_1\},~g_1\in
G_1$. For convenience we will assume that $K$ contains the
compact set $U_1\subset G_1$ defined in the proof of
Proposition \ref{prop:C<infty} below.

We define the (leafwise $G_1$--) energy of a $\gC$--equivariant function $\varphi:G\to
X$ to be
$$
 E(\varphi )=\int_{\gO\times G_1}h(g_1)d\big(\varphi (\go),\varphi (\go g_1)\big)^2=
 \int_{(\gC\backslash G)\times G_1}h(g_1)d\big(\varphi (g),\varphi (gg_1)\big)^2.
$$
Note that the energy $E$ is convex and $G_2$--invariant from the right, i.e.
$E(\varphi)=E(\varphi (\cdot g_2))$ for any $g_2\in G_2$. It is
also easy to check that if $\varphi\in L_2(\gO,X)$ then $E(\varphi
)<\infty$ and that $E$ is continuous on $L_2(\gO,X)$. Let
$M=\inf\{ E(\varphi ):\varphi\in L_2(\gO,X)\}$.

\begin{defn}
$\varphi\in L_2(\gO,X)$ is called harmonic if $E(\varphi )=M$.
\end{defn}

\begin{thm}\label{thm:existence}
There exists a harmonic map.
\end{thm}

If $\gC$ has a global fixed point $y_0$ then the constant map
$\gph (g)\equiv y_0$ is $\gC$--equivariant with energy $0$, and
hence harmonic. For the rest of this section, we will assume that
$\gC$ has no global fixed point in $X$.

Let us fix the point $x_0\in X$ and denote by $x_0$ also the element in $L_2(\gO,X)$ which sends $\gO$ to $x_0$. For $\gph\in L_2(\gO,X)$ let 
$$
 \|\varphi\|:=\rho (\varphi,x_0)=[\int_\gO d(\varphi (\go
 ),x_0)^2]^{\frac{1}{2}}.
$$

For each $n$, let $\varphi_n\in L_2(\gO,X)$ be a map satisfying:
\begin{itemize}
\item $E(\varphi_n)\leq M+\frac{1}{n}$, and

\item $\|\varphi_n\|\leq\inf\{\|\varphi\| :\varphi\in
L_2(\gO,X),E(\varphi )\leq M+\frac{1}{n}\}+\frac{1}{n}$.
\end{itemize}

\begin{prop}\label{prop:C<infty}
The maps $\gph_n$ are uniformly bounded, i.e.
$\sup\|\varphi_n\|<\infty$.
\end{prop}

\begin{proof}
Let
$$
 \mathcal{F}_n=\inf\{\|\varphi\| :\varphi\in
 L_2(\gO,X),E(\varphi )\leq M+\frac{1}{n}\}.
$$
Note that
$\mathcal{F}_n\leq\|\varphi_n\|\leq\mathcal{F}_n+\frac{1}{n}$.

Let $\ti\gO =\gS\cdot\gO,~U=\gO^{-1}\ti\gO$, and
$U_i=\pi_i(U),~i=1,2$. For $u\in U$ let $u_i=\pi_i(u)\in U_i$,
where $\pi_i:G\to G_i$ is the canonical projection. We may
normalize the Haar measure $\mu_G$ of $G$ so that $\mu_G(\gO )=1$.

Let $b$ be the constant from Lemma \ref{lem:b>0} applied to the
function $d_\gS$, then
$$
 \frac{(b\|\varphi_n\|)^2}{4}\leq \frac{1}{4}\int_\gO d_\gS(\varphi_n(\omega))^2\leq
 \inf_{y\in X}\int_\gO\max_{\gc\in\gS}d(y,\gc\cdot\gph_n(\omega))^2 \leq \inf_{y\in
 X}\int_\gO\sum_{\gc\in\gS}d(y,\gc\cdot\gph_n(\omega))^2
$$
where the second $\leq$ follows from the triangle inequality since $\gS$ contains 1.
By the definition of $\ti\gO$, the last term is equal to
$$
 \inf_{y\in X}\int_{\ti\gO}d(y,\gph_n(\ti\omega))^2\leq
 \int_{\gO\times\ti\gO}d(\varphi_n(\go ),\varphi_n(\ti\go ))^2\leq
 \int_{\gO\times U}d(\varphi_n(\go ),\varphi_n(\go u))^2
$$
here the first $\leq$ holds since $\mu_G(\gO)=1$, and the second since $\ti\gO\subset wU,~\forall w\in\gO$.
By the triangle inequality, the last term is bounded by
$$
 \int_{\gO\times U}[d(\varphi_n(\go ),\varphi_n(\go u_2))+d(\varphi_n(\go u_2),\varphi_n(\go u_2u_1))]^2
$$
$$
 \leq
 2\int_{\gO\times U}d(\varphi_n(\go ),\varphi_n(\go u_2))^2+
 2\int_{\gO\times U_2\times U_1}d(\varphi_n(\go u_2),\varphi_n(\go u_2u_1))^2.
$$
Now the second summand is bounded above by
$2\mu_2(U_2)E(\varphi_n)\leq 2\mu_2(U_2)(M+1)$, since $h|_{U_1}\geq 1$.
Hence, if we assume that $\frac{b^2}{4}\|\varphi_n\|^2\geq
4\mu_2(U_2)(M+1)$ then we have
$$
 \int_{\gO\times U}
 d(\varphi_n(\go ),\varphi_n(\go u_2))^2\geq
 \frac{(b\|\varphi_n\|)^2}{16}
$$
and therefore, for some $u_2\in U_2$
$$
 \rho (\varphi_n,\varphi_n(\cdot u_2))^2=\int_\gO
 d(\varphi_n(\go),\varphi_n(\go u_2))^2
 \geq\frac{1}{\mu (U)}\frac{b^2}{16}\|\varphi_n\|^2.
$$
On the other hand the relatively compact set $\gO U_2$ is
contained in $\ti\gS\cdot\gO$ for some finite set
$\ti\gS\subset\gC$, and we have $\|\varphi_n(\cdot u_2)\|\leq
\|\varphi_n\|+d_{\ti\gS}(x_0)$, and if we assume further that
$\|\gph_n\|\geq d_{\ti\gS}(x_0)$ and take
$\gep_0:=\frac{b}{8\sqrt{\mu (U)}}$ then we also have $\rho
(\varphi_n,\varphi_n(\cdot u_2))\geq\gep_0(\|\gph_n\|+d_{\ti\gS}(x_0))$.

Let $\varphi_n'=\frac{\varphi_n+\varphi_n(\cdot u_2)}{2}$ then
$$
 E(\varphi_n')\leq E(\varphi_n)\leq M+\frac{1}{n}
$$
since $E$ is convex and $G_2$--invariant from the right. Besides, if we let $\gt(\gep)>0$ be a linear lower bound for the modulus of convexity function of $L_2(\gO,X)$, as in the definition of uniform convexity (see Proposition \ref{prop:Lp}), then 
$$
 \|\varphi_n'\|\leq (\|\varphi_n\|+d_{\ti\gS}(x_0))-\gt(\gep_0)
 (\|\varphi_n\|+d_{\ti\gS}(x_0))\leq (\mathcal{F}_n+1+d_{\ti\gS}(x_0))-
 \gt(\gep_0)(\|\varphi_n\|+d_{\ti\gS}(x_0)),
$$
and since $\mathcal{F}_n\leq\|\varphi_n'\|$, this implies that
$\gt(\gep_0)(\|\varphi_n\|+d_{\ti\gS}(x_0))\leq
d_{\ti\gS}(x_0)+1$. Thus $\|\varphi_n\|$ is bounded independently of $n$.
\end{proof}

We are now able to finish the proof of Theorem
\ref{thm:existence}. Note that we may assume that the positive function $\gt(\gep)$ is monotonic (non-decreasing) in $\gep$. Let $\mathcal{F}=\lim\mathcal{F}_n$ and let $m>n$, then
$$
 \mathcal{F}_n\leq \|\frac{\varphi_n+\varphi_m}{2}\|\leq
 (\mathcal{F}+1/n)\big(1-\gt(\frac{\rho
 (\varphi_n,\varphi_m)}{\mathcal{F}+1/n})\big)
$$
therefore
$\gt(\frac{\rho (\varphi_n,\varphi_m)}{\mathcal{F}+1/n})\to 0$, 
which force 
$\rho(\gph_n,\gph_m)\to 0$ when $n,m\to\infty$.
Thus $\varphi_n$ is a Cauchy sequence.
Since $L_2(\gO,X)$ is complete we can take $\varphi
:=\lim\varphi_n$ and by continuity of $E$, $\varphi$ is a harmonic
map. $\qed$



\section{The proof of Theorem \ref{thm1}}\label{proof}
\subsection{The case $n=2$}\label{subsec:n=2}
Let $G,\gC,X$ be as in Theorem \ref{thm1}, and suppose that we are
given a reduced $\gC$--action on $X$ without global fixed points.
We will first deal with the case $n=2$, i.e. $G=G_1\times G_2$.

Observe the following facts about our harmonic maps:
\begin{itemize}
\item If $\gph$ is a harmonic map then so is $\gph(\cdot g_2)$ for
any $g_2\in G_2$, because the energy $E$ is $G_2$--invariant from
the right.

\item If $\gph$ and $\psi$ are harmonic maps then so is
$\frac{\gph+\psi}{2}$ (by convexity of $E$).

\item If $\gph$ and $\psi$ are harmonic maps then $[\gph (g),\gph
(gg_1)]~\|~[\psi (g),\psi (gg_1)],\forall g\in G,g_1\in G_1$ (by BNPC of $X$).
\end{itemize}

Let us now fix a harmonic map $\gph$. From the first and the third
facts we see that
$$
 [\gph (g),\gph(gg_1)]~\|~[\gph (gg_2),\gph(gg_1g_2)],
$$
and by Lemma \ref{lem:parallelogram}
$$
 [\gph (g),\gph(gg_2)]~\|~[\gph (gg_1),\gph(gg_1g_2)],
$$
for almost all $g\in G,g_1\in G_1,g_2\in G_2$.

\begin{lem}\label{lem:const-dist}
For any fixed $g_1\in G_1$ the function $g\mapsto d(\gph (g),\gph
(gg_1))$ is essentially constant. Similarly, for any fixed $g_2\in
G_2$ the function $g\mapsto d(\gph (g),\gph (gg_2))$ is
essentially constant.
\end{lem}

Lemma \ref{lem:const-dist} is similar to a statement implanted in the proof of the main theorem in \cite{Monod}.

\begin{proof}
The first (resp. second) function is measurable, $\gC$--invariant
from the left and $G_2$ (resp. $G_1$) invariant from the right.
Since $\gC$ is irreducible, $G_2$ (resp. $G_1$) acts ergodically
on $\gC\backslash G$ from the right. The result follows.
\end{proof}

\begin{cor}\label{cor:continuity}
The harmonic map $\gph$ is essentially continuous.
\end{cor}

\begin{proof}
It follows from Lemma \ref{lem:const-dist} that
$$
d(\gph (g),\gph
(gg_1))=\rho (\gph,\gph (\cdot g_1))\textrm{ and }d(\gph (g),\gph
(gg_2))=\rho (\gph,\gph (\cdot g_2)),
$$
for almost any $g\in
G,g_i\in G_i$. Since the (right) action of $G$ on
$L_2(\gC\backslash G,X)$ is continuous\footnote{This fact is a
well known when $X$ is replaced by $\BR$, and is easily seen to be
true for any metric space $X$.}, $\gph$ is essentially uniformly
(on $G$) continuous from the right along $G_1$ and $G_2$, and
hence along $G$. It follows that $\gph$ is essentially continuous.
\end{proof}

By changing $\gph$ on a set of measure $0$ we can assume that it
is actually continuous. All harmonic maps considered bellow will
be assumed to be continuous rather than essentially continuous.

For the proof of Theorem \ref{thm1} we will distinguish between
two cases:
\begin{enumerate}
\item $\gph (G_2)$ is bounded.

\item $\gph (G_2)$ is unbounded.
\end{enumerate}

In case $(1)$ we will show that there is a $G_2$--invariant
harmonic map $\gph_0$ (perhaps different from $\gph$). In case
$(2)$, we will show that $\gph$ is $G_1$--invariant. We will then
conclude that in case $(i)$ (where $i=1$ or $2$) the $\gC$ action extends to a
continuous $G$--action which factors through $G_i$.

Before we start let us note that in a uniquely geodesic metric space
$X'$, the closed convex hull $\overline{\text{conv}(Y)}$ of a
subset $Y\subset X'$, which by definition is the minimal closed
convex subset of $X'$ containing $Y$ (which is also the
intersection of all such sets) can be constructed recursively as
follows: Define $Y_0=Y$, and $Y_n=\cup\{\frac{x+y}{2}:x,y\in
Y_{n-1}\}$. Then $Y_n\supset Y_{n-1}$ because we can take $x=y$,
and $\overline{\text{conv}(Y)}=\overline{\cup_{n=0}^\infty Y_n}$.

\medskip

$(1)$ If $\gph (G_2)$ is bounded, then, as follows from Lemma
\ref{lem:const-dist}, the set
$$
 G_2\cdot\gph=\{ \gph (\cdot g_2):g_2\in G_2\}
$$
is bounded in $L_2(\gO,X)$. By the constructive
description of $\overline{\text{conv}(G_2\cdot\gph )}$ we see that
it consists of harmonic maps. Since
$\overline{\text{conv}(G_2\cdot\gph )}$ is bounded and convex and
since $L_2(\gO,X)$ is uniformly convex, there is a unique relative circumcenter, i.e. a unique point $\gph_0\in\overline{\text{conv}(G_2\cdot\gph )}$ which minimize
$$
 \sup\{\rho (\gph_0,\gph'):\gph'\in\overline{\text{conv}(G_2\cdot\gph
 )}\}.
$$
(The existence and uniqueness of a relative circumcenter of a closed bounded convex subset of a complete WUC space follow from the existence and uniqueness of the usual circumcenter by ignoring the complement of the set.)

Since $\overline{\text{conv}(G_2\cdot\gph )}$ is $G_2$--invariant,
the function $\gph_0$ is $G_2$--invariant.

\medskip

$(2)$ Suppose that $\gph (G_2)$ is unbounded, and let
$Y=\overline{\text{conv}(\gph (G_2))}$. Then the Hausdorff
distance $\text{Hd}(\gc\cdot
Y,Y)<\infty$ for any $\gc\in\gC$. Indeed $\gc\cdot\gph (g_2)=\gph
(\gc_2g_2\gc_1)$ where $\gc =(\gc_1,\gc_2)$. Hence $d(\gc\cdot\gph
(g_2),\gph (\gc_2g_2))$ is a constant depending on $\gc_1$ as
follows from Lemma \ref{lem:const-dist}. It follows, since the
action $\gC\circlearrowleft X$ is assumed to be reduced, that
$Y=X$.

By the constructive description of the closed convex hull we see
that $X=Y=\overline{\cup_{n=0}^\infty Y_n}$ where each point of
$\cup_{n=0}^\infty Y_n$ is of the form $\psi (1)$ for some
harmonic function $\psi$ in the convex hull of $G_2\cdot\gph$. Let
$g_1\in G_1$. Since the segments $[\psi (1),\psi
(g_1)],~\psi\in\text{conv}(G_2\cdot\gph )$ are all parallel to
each other, the map $\psi (1)\mapsto\psi (g_1)$ extends to a
Clifford isometry on $\overline{\{\psi
(1):\psi\in\text{conv}(G_2\cdot\gph )\}}$, namely on $X$. Since
$X$ has no non-trivial Clifford isometries, we get that $\gph
(g_2)=\gph (g_1g_2)$ for any $g_2\in G_2$. Since $g_1$ is
arbitrary, we get that $\gph$ is $G_1$--invariant.

\medskip

We showed that there is a harmonic map $\gph_0$ which is either
$G_1$ or $G_2$ invariant. Suppose it is $G_2$--invariant. Then,
since $\gph_0$ is $\gC$--equivariant and continuous, the orbit map
$\gc\mapsto\gc\cdot x$ is continuous for any $x\in\gph_0(G)$,
where the topology on $\gC$ is the (not necessarily Hausdorff) one
induced from $G_1$ (equivalently, we can consider $\gC$ as a dense
subgroup of $G_1$). It follows that the set
$$
 \{ x\in X:\text{the obit map~}~\gc\mapsto\gc\cdot x~\text{~is
 continuous with respect to
 the~}~G_1\text{--topology}\}
$$
is nonempty. Since that set is also closed convex and
$\gC$--invariant, it follows from $C$--minimality that it is the
whole space $X$. Thus the orbit map is continuous for any point
$x\in X$, when $\gC$ is considered with the topology induced from
$G_1$, so we can define the action of $G$ on $X$ by
$$
 g\cdot x:=\lim_{\pi_1(\gc )\to\pi_1 (g)}\gc\cdot x.
$$


\subsection{The case $n>2$}\label{n>2}

Let now $n\geq 2$ be general, $G=\prod_{i=1}^nG_i$, and define $n-1$ energies $E_1,E_{1,2},\ldots,E_{1,\ldots,n-1}$ as follow:
$$
 E_{1,\ldots,k}(\gph):=\int_{(\gC\backslash G)\times G_1\times\ldots\times G_k}h_1(g_1)\ldots h_k(g_k)d(\gph (g),\gph (gg_1\ldots
 g_k))^2,
$$
where $h_i(g_i)=e^{-|g_i|^2}$, $|~|$ being the norm with respect
to the word metric on $G_i$ associated with some compact
generating set.

The set $H_1$ of $E_1$--harmonic maps is closed convex, and by
Theorem \ref{thm:existence} non-empty. Furthermore, there exist
some $\gph\in H_1$ which minimize $E_{1,2}$. To see this, argue as
in the proof of Theorem \ref{thm:existence}, taking the $\gph_n$
to be $E_1$--harmonic and letting $\prod_{i>2}G_i$ play the rule of $G_2$ in
\ref{thm:existence}. Call such a function $E_{1,2}$--harmonic. More generally, for all $k<n$
call a $\gC$--equivariant map $E_{1,\ldots,k}$--harmonic if it minimizes $E_{1,\ldots,k}$ among the $E_{1,\ldots,k-1}$ harmonic maps. By repeating the argument above finitely many times, one proves the existence of an $E_{1,\ldots,n-1}$--harmonic map $\gph$.
Since $\gph (G)$ is $\gC$--invariant it must be unbounded. Let $1\leq k\leq n$ be the smallest integer such that $\gph (G_{k+1}\times\ldots\times G_n)$ is bounded. Then $\gph (G_k)$ is unbounded, and by replacing $\gph$ by the relative circumcenter of
$$
 \overline{\text{conv}\{\gph (\cdot g):g\in G_{k+1}\times\ldots\times G_n\}}
$$
one gets an $E_{1,\ldots,k'}$--harmonic map $\gph_0$ which is $G_i$ invariant for all $i>k$, where $k'=\min\{k,n-1\}$. As in case (2) of the proof of \ref{thm:existence} one deduces, since $\gph_0$ is also $E_{1,\ldots,k-1}$--harmonic (in case $k>1$), that $\gph$ is also $G_i$--invariant for all $i<k$. Therefore one can use $\gph_0$ to extend the $\gC$--action to a continuous $G$ action which factors through $G_k$.
\qed


\section{The proof of Theorem \ref{thm2}} For simplicity assume
again that $n=2$. If $d_\gS\nrightarrow\infty$ where $\gS$ is a
finite generating set, then $\gC$ has a fixed point in $\partial
X$. Assume therefore that $d_\gS\to\infty$. By Lemma
\ref{lem:C-min} there is an unbounded closed convex
$\gC$--invariant subset on which the action is $C$--minimal.
Replacing $X$ with such a subset we may assume that the
action on $X$ is $C$--minimal. Moreover, since $d_\gS\to\infty$ we
have a harmonic map $\gph$, and we can argue as in paragraph
\ref{subsec:n=2}. In case $(1)$, when $\gph (G_2)$ is bounded, we
obtain that the $\gC$--action extends to a continuous $G$--action
on $X$ which factors through $G_1$. In case $(2)$, when $\gph
(G_2)$ is unbounded, then for any $g_1\in G_1$ the map $\gph
(g_2)\mapsto\gph (g_1g_2)$ extends to a parallel translation from
$\overline{\text{conv}(\gph (G_2))}$ onto
$\overline{\text{conv}(\gph (g_1G_2))}$. Hence, in that case $\gC$
preserves $\partial\overline{\text{conv}(\gph (G_2))}$ and the
action on it factors through $G_2$. Finally, since $\gph$ is
continuous and $X$ is proper, the $\gC$--action on
$\partial\overline{\text{conv}(\gph (G_2))}$ is continuous with
respect to the $G_2$--topology, and extends to a continuous
$G_2$--action.

When the number of factors $n$ is greater than two, one argues as in Paragraph \ref{n>2}.
\qed


\section{The proof of Theorem \ref{thm3}}
By Lemma \ref{S-min}, up to replacing $X$ with a closed non-empty
geodesically complete subset, we may assume that $X$ is a minimal
complete uniquely geodesically complete for the $\gC$--action.

For a subset $A\subset X$ of cardinality $\geq 2$
let span$(A)$ denote the minimal set containing $A$ with the
geodesic extension property. One can construct span$(A)$
recursively by defining $A_0=A$ and $A_{m+1}=\cup\{\overline{xy}:x,y\in
A_m\}$, where $\overline{xy}$ is the geodesic
containing $x,y$, and taking the union span$(A)=\cup A_m$.

For the sake of simplicity let us assume again that $n=2$, i.e.
$G=G_1\times G_2$. One can extend the argument below to any $n\geq 2$ using the energies $E_{1,\ldots,k},~k=1,\ldots n-1$ as in Paragraph \ref{n>2}.

Let $\gph :G\to X$ be a harmonic map. If $\gph
(G_i)$ is a single point, for $i=1$ or $2$, then $\gph$ is
$G_i$--invariant and the action extends to a $G$--action which
factors through $G_{3-i}$. Thus we may assume that $|\gph (G_i)|\geq 2$,
for $i=1,2$.

Let $H\leq L_2(\gO,X)$ be the subset of harmonic maps.

\begin{lem}\label{lem:H-gc}
$L_2(\gO,X)$ is complete CAT(0) and uniquely geodesically
complete, and so is its subset $H$.
\end{lem}

\begin{proof}
The first statement is straightforward. The second one follows,
using Lemma \ref{lem:parallelogram}, from the fact that a
$\gC$--equivariant map $\psi :G\to X$ is harmonic if and only if for almost
all $g\in G,g_1\in G_1$ the segment $[\psi (g),\psi (gg_1)]$ is
parallel to $[\gph (g),\gph (gg_1)]$.
\end{proof}

Let
$$
 X_1=\overline{\text{span}(\gph (G_1))}~\text{~and~}~X_2=\{\psi
 (1):\psi\in H\}.
$$
By definition $X_1$ is closed and geodesically complete, and it
follows from Lemma \ref{lem:H-gc} that also $X_2$ is. We will show
that $X=X_1\times X_2$. For this we need:

\begin{lem}\label{lem:spliting}
Let $\gph_1,\gph_2$ be two harmonic maps, then
$$
 \text{conv}\big(\overline{\text{span}(\gph_1(G_1))}\cup
 \overline{\text{span}(\gph_2(G_1))}\big)
 \cong \overline{\text{span}(\gph (G_1))}\times
 [0,d(\text{span}(\gph_1(G_1)),\text{span}(\gph_2(G_1))\big) ] .
$$
\end{lem}

\begin{proof}[Proof of Lemma \ref{lem:spliting}]
Since $X$ is CAT(0), parallelity is a transitive relation on the set of two sided infinite
geodesics in $X$ (cf. \cite{Br-Ha} p. 183), and we can use this property
to extend the map $\gph_1 (g_1)\mapsto\gph_2(g_1)$ to a parallel
isometry
$T:\overline{\text{span}(\gph_1(G_1))}\to\overline{\text{span}(\gph_2(G_1))}$,
parallel means that 
$$
[x,T(x)]~\|~[y,T(y)],~\forall
x,y\in\overline{\text{span}(\gph_1(G_1))}.
$$

The lemma would follow from the flat strip theorem (see
\cite{Br-Ha}, p. 183; The flat strip theorem is stated for parallel geodesics, but the proof extends straightforwardly to parallel geodesically complete sets) once we show that
$d(x,\overline{\text{span}(\gph_2(G_1))})$ is constant over
$\overline{\text{span}(\gph_1(G_1))}$. Suppose in contrary that
there are $x,y\in\overline{\text{span}(\gph_1(G_1))}$ with
$$
 d(x,\overline{\text{span}(\gph_2(G_1))})<d(y,\overline{\text{span}(\gph_2(G_1))}).
$$
Let $c(t)$ be the geodesic with $c(0)=x,c(d(x,y))=y$. Since the
function
$$
t\mapsto d\big(
c(t),\overline{\text{span}(\gph_2(G_1))}\big)
$$
is convex on $\BR$,
$d\big( c(t),\overline{\text{span}(\gph_2(G_1))}\big)\to\infty$
when $t\to +\infty$, contradicting the fact that $d(c(t),T(c(t)))$
is constant.
\end{proof}

Let 
$$
Y=\overline{\text{conv}\{\bigcup \psi (G):\psi\in H\}}=\overline{\text{conv}\{\bigcup \psi (G_1):\psi\in H\}}.$$ 
Using again the fact that parallelity is a transitive relation and
the flat strip theorem, one sees that every point in $Y$ is
contained in an isometric parallel copy of $X_1$. It follows that
$X_1$ is a direct factor of $Y$ (In case $X_1$ is a single geodesic, this is the content of \cite{Br-Ha} II.2.14 which again extends straightforwardly to the more general case where $X_1$ is geodesically complete). Furthermore, since $Y$ is
$\gC$--invariant it follows that $X=\overline{\text{span}(Y)}$,
which allows one to show that every point in $X$ is contained in
an isometric parallel copy of $X_1$. It follows that $X_1$ is a
direct factor of $X$.

Next we claim that for any $\gph_1,\gph_2\in H$ the projection
$$
p:\overline{\text{span}(\gph_1(G_1))}\to\overline{\text{span}(\gph_2(G_1))}
$$
extending the parallel translation $\gph_1(g_1)\mapsto\gph_2(g_1)$
is an orthogonal projection. Indeed, if that was not the case then
the orthogonal projection
$\pi:\overline{\text{span}(\gph_2(G_1))}\to\overline{\text{span}(\gph_1(G_1))}$
composed with $p$ would be a non-trivial Clifford isometry on
$\overline{\text{span}(\gph_1(G_1))}$. To see this note that for any
$x_1,x_2\in \overline{\text{span}(\gph_1(G_1))}$
$$
 [x_1,x_2]~\|~[p(x_1),p(x_2)]~\|~[\pi (p(x_1)),\pi (p(x_2))]
$$
which, by transitivity (since $X$ is uniquely geodesically complete parallelity is a transitive relation on segments as well) and Lemma \ref{lem:parallelogram}, implies
that 
$[x_1,\pi (p(x_1))]~\|~[x_2,\pi (p(x_2))]$.  
This however
would imply that $\overline{\text{span}(\gph_1(G_1))}$ has a
euclidian factor. Now since
$\overline{\text{span}(\gph_1(G_1))}$ is isometric to $X_1$ and
since $X_1$ is a direct factor of $X$, this contradicts our
assumptions on $X$.

It follows in particular that the set 
$Z=\overline{\{\bigcup\text{span}(\psi (G_1)):\psi\in H\}}$ 
decomposes as a direct product 
$$
 Z=\text{span}(\gph (G_1))\times\{\psi (1):\psi\in H\}=X_1\times X_2
$$ 
and hence it is geodesically complete. Since $Z$ is also $\gC$--invariant, we derive from minimality that $X=Z$ and hence $X=X_1\times X_2$.

Finally we claim that $G_i$ acts by isometries on $X_i$, for $i=1,2$.
First note that if $\psi\in H$ and $g_2\in G_2$ then also $\psi
(\cdot g_2)\in H$ and since $X_2=\{\psi (1):\psi\in H\}$ we can
define an action of $G_2$ on $X_2$ by
$$
 g_2\cdot\psi (1)=\psi (g_2).
$$
We need to show that this action is well defined and is by
isometries. These two claims follows from the fact that $\gC$ is
irreducible in $G$. Indeed, any $g_2\in G_2$ can be approximated by
$\gc g_1$ where $\gc\in\gC,g_1\in G_1$, and we know that for any
two harmonic maps $\gph_1,\gph_2$ the segments
$[\gph_1(1),\gph_1(g_1)]$ and $[\gph_2(1),\gph_2(g_1)]$ are
parallel. Since $\gc$ is an isometry and $\gph_i$ are $\gC$--equivariant, we derive that 
$$
 d(\gph_1(\gc g_1),\gph_2(\gc g_2))=d(\gph_1(g_1),\gph_2(g_1))=d(\gph_1(1),\gph_2(1)),
$$ 
and by passing to limit, $d(\gph_1(g_2),\gph_2(g_2))=d(\gph_1(1),\gph_2(1))$.

Similarity one shows that $G_1$ acts on $X_1$ by isometries (with
$g_1\cdot\gph (g_1')=\gph (g_1g_1')$) by approximating $g_1$ with
$\gc g_2\in\gC G_2$.

This produces an action of $G$ on $X=X_1\times X_2$ which extends
the $\gC$--action and is clearly continuous. \qed


\section{Non-uniform weakly cocompact $p$--integrable lattices}\label{sec:non-uniform}

Let $G$ be a locally compact group, $\gC$ a lattice in $G$ and
$\gO$ a right fundamental domain. We define a map $\chi :G\to\gC$
by the rule $g\in\chi (g)\gO$. Suppose that $\gC$ is generated by
a finite set $\gS$ and let $|~|:\gC\to\BN$ be the word norm
associated to $\gS$. Let $p>1$. In analogy to \cite{Shalom}, we
will say that $\gC$ is $p$--integrable if for any element $g\in G$
the function $\go\mapsto |\chi (\go g)|$ belongs to $L_p(\gO )$.
This assumption ensure that whenever $\gC$ acts by isometries on a
metric space $X$, the space $L_p(\gO,X)$ is invariant under $G$, i.e. if $\gph :G\to X$ is
$\gC$--equivariant and $\gph|_\gO\in L_p(\gO,X)$ then also $\gph
(\cdot g)\in L_p(\gO,X)$ for any $g\in G$.  Note that the property of being $p$--integrable is
independent of the generating set $\gS$, however, it does depend
on the choice of $\gO$. When this condition is satisfied we shall say that $\gO$ is 
$p$--{\it admissible}.

Let $L_p^0(\gC\backslash G)$ denotes the codimension one subspace
of $L_p(\gC\backslash G)$ of function with $0$ mean. Following
[\cite{Margulisb}, III.1.8] we will say that $\gC$ is weakly
cocompact if the right regular representation of $G$ on
$L_p^0(\gC\backslash G)$ does not almost have invariant vectors. This is equivalent to each of the following:
\begin{enumerate}
\item If $f_n\in L_p(\gC\backslash G)$ are normalized
asymptotically invariant positive functions then $(f_n)$ converges to a
constant function.

\item If $f_n\in L_p(\gC\backslash G)$ are normalized positive functions
such that for any compact $K\subset G$, $\int_K\|f_n-f_n(\cdot
k)\|_p\to 0$ then $(f_n)$ converges to a constant function.
\end{enumerate}
Moreover, using the Mazur map $M_{p,q}:L_p(\gC\backslash G)\ni f\mapsto |f|^{\frac{p}{q}}\text{sign}(f)\in L_q(\gC\backslash G)$ which intertwines the $G$ actions and is uniformly continuous (c.f. \cite{BL} Theorem 9.1) one can show that this property is independent of $1<p<\infty$.

The following extends the superrigidity results from the previous sections to non-uniform
lattices which are $p$--integrable and weakly cocompact. 
This is analogous to [\cite{Monod}, Theorem 7] which gives a similar statement for 
actions of weakly cocompact 2--integrable non-uniform lattices on CAT(0) spaces.

\begin{thm}\label{thm:non-uniform}
Theorems \ref{thm1}, \ref{thm2}, \ref{thm3} and \ref{thm:CL} remain true for
finitely generated non-uniform lattices $\gC$ provided they are
weakly cocompact and $p$--integrable for some $1<p<\infty$.
\end{thm}

For the sake of simplicity, let us assume that $G=G_1\times G_2$
is a product of two factors. Fix $1<p<\infty$ such that $\gC$ is
$p$--integrable with respect to $\gO$, and normalize the Haar measure so that $\mu (\gO )=1$. Then
$G$ acts measurably  and hence continuously from the right on $L_p(\gO,X)$. It also follows from
$p$--integrability that one can choose a measurable function $h:G_1\to\BR_{>0}$ such that the energy
$E(\gph):=\int_{\gO\times G_1}h(g_1)d(\gph (g),\gph (gg_1))^p$ is
finite for every $\gph\in L_p(\gO,X)$. Since $G$ acts continuously
on $L_p(\gO,X)$ we can take $h$ with $\inf_{k\in
K_1}h(k)>0$ for every compact $K_1\subset G_1$.

\medskip

The only place in the previous sections that
compactness of $\gO$ was used is the proof of Proposition
\ref{prop:C<infty}. Hence we should only justify why the function
$\gph_n$, chosen as in Section \ref{sec:existence}, are uniformly
bounded when $\gO$ is not relatively compact but $\gC$ is
$p$--integrable and weakly cocompact. Let $\rho$ denote the
distance on $L_p(\gO,X)$, $\gph_{x_0}$ the $\gC$--equivariant
function sending $\gO$ to $x_0\in X$, and for $\gph\in L_p(\gO,X)$
set $\|\gph\|=\rho (\gph,\gph_{x_0})$.

\begin{lem}\label{lem:K,beta}
There exists a compact subset $K\subset G$ and a positive constant
$\gb>0$ such that any $\gph\in L_p(\gO,X)$ satisfies
$\int_{K}\rho(\gph,k\cdot\gph)>\gb\|\gph\|-\frac{1}{\gb}$.
\end{lem}

\begin{proof}
Assuming the contrary, as $G$ is $\gs$--compact, one can find a
sequence $\psi_n\in L_p(\gO,X)$ with $\|\psi_n\|\to\infty$ and
$\frac{1}{\|\psi_n\|}\int_{K'}\rho(\psi_n,k\cdot\psi_n)\to
0$ for every compact $K'$. Let
$f_n'(g)=d(\psi_n(g),\gph_{x_0}(g))$ and
$f_n=\frac{f_n'}{\|f_n'\|_p}$. It is straightforward to verify  that these $(f_n)$ satisfy the
condition (2) above and hence $(f_n)$ converges to a constant function.
Let $K_0\subset\gO$ be a compact subset with positive measure. Then
if $n$ is sufficiently large, for every $k$ in some subset
$K_n\subset K_0$ of at least half the measure of $K_0$,
$d(x_0,\psi_n(k))>\frac{1}{2}\|\psi_n\|$. Taking
$K=K_0^{-1}\gS\cdot K_0$ and using Lemma \ref{lem:b>0} one gets a
constant $\gb>0$ for which
$\int_{K}\rho(\psi_n,k\cdot\psi_n)>\gb\|\psi_n\|$
for all sufficiently large $n$, in contrary to the assumptions on
$\psi_n$.
\end{proof}

Consider now the functions $\gph_n$ defined as in Section
\ref{sec:existence}. Replacing $K$ by a larger compact set, we may
assume it is of the form $K=K_1\times K_2$ with $K_i\subset G_i$.
As in Section \ref{sec:existence} the energy of $\gph_n$ bounds
its variation along the $G_1$ factor, so we conclude that for some
positive constant $\gb'$
$$
 \int_{K_2}\rho(\gph_n,k_2\cdot\gph_n)>\gb'\|\gph_n\|-\frac{1}{\gb'},
$$
and hence for some $k_2\in K_2$ and another positive constant $\gep'$ we have
$\rho (\gph_n,\gph_n(\cdot k_2))> \gep'\|\gph_n\|-\frac{1}{\gep'}$. Finally since $k_2$ belongs to the compact set $K_2$ the norm of $\gph_n(\cdot k_2)$ is bounded by the norm of $\gph_n+$ some constant. Define $\gph_n'=\frac{\gph_n+\gph_n(\cdot k_2)}{2}$ then $E(\gph_n')\leq E(\gph_n)$ since $E$ is $G_2$--invariant and convex, and $\|\gph_n\|-\|\gph_n'\|\geq \gd'\|\gph_n\|-\frac{1}{\gd'}$ for some positive constant $\gd'$, by uniform convexity of $L_p(\gO,X)$. This, together with the second property of $\gph_n$ implies that its norm must be bounded independently of $n$.

\begin{rem}\label{rem:p-int}
(i) It is possible to show that if $G$ is a semisimple real
Lie group without compact factors not locally isomorphic to
$SL(2,\mathbb{R})$, then any lattice in $G$ is $p$--integrable for some $p>1$. Similarly, all irreducible lattice in higher rank groups over local fields are $2$--integrable (c.f. \cite{Shalom}). R\'emy \cite{Remy} showed that all Kac--Moody lattices are $2$--integrable.

(ii) Every lattice in a semisimple Lie group over a local field is weakly cocompact [\cite{Margulisb}, III.1.12]. Clearly, if $G$ has Kazhdan property (T) then any lattice in $G$ is weakly cocompact.
\end{rem}


\section{Superrigidity for commensurability subgroups}\label{sec:commensurators}

Let $G$ be a locally compact, compactly generated group and $\Gamma$ a
cocompact lattice in $G$. Let $\Lambda$ be a subgroup of
\[
Comm_{G}(\Gamma):=\{g\in G:g\Gamma g^{-1}\text{ and }\Gamma\text{ are
commensurable}\}
\]
containing $\Gamma$ which is dense in $G$. Let $X$ be a complete uniformly convex BNPC metric space. Assume that $\Lambda$ acts by isometries on $X$ such that any
subgroup $\Gamma_{0}\leq\gD$ commensurable to $\Gamma$ satisfies
$d_{\gS_0}\to\infty$ where $\gS_0$ is a finite generating set of
$\gC_0$, and has no parallel orbits\footnote{Note that if $X$ is uniquely geodesically complete then the assumption that there are no parallel orbits follows from the assumption that $d_{\gS_0}\to\infty$.\label{fn:disp-->para}}, i.e. for any two distinct points $x,y\in X$ there is $\gc\in\gC_0$ (equivalently, there is $\gc\in \gS_0$) for which $[\gc\cdot x,\gc\cdot
y]\nparallel [x,y]$. Assume moreover that the action is C--minimal. The following superrigidity theorem was
proved in \cite{Margulis} (under the weaker assumption that $X$ is WUC rather than UC):

\begin{thm}\label{thmcommensurator}
Under the above assumptions, the $\Lambda$--action extends uniquely to a
continuous isometric $G$--action.
\end{thm}

For a subgroup $\Gamma_{0}\leq\gD$ commensurable to $\gC$ we
define a $\Gamma_{0}$\emph{--harmonic map} to be a
$\Gamma_{0}$--equivariant map in $L_{2}(\Omega_0,X)\ $ which
minimizes
$$
I_{\Gamma_{0}}(\varphi)=\frac{1}{\mu (\gO_0)}\int_{\Omega_0\times
G}h(\go^{-1}g)d(\varphi(\go),\varphi(g))^{2}
$$
$$
=\frac{1}{\mu
(\gO_0)}\int_{\gC_0\backslash (G\times G)}h(g_1^{-1}g_2)d(\gph
(g_1),\gph (g_2))^2,
$$
where $\Omega_0$ is a relatively compact fundamental domain
for $\Gamma_{0}$ in $G$ and $h$ is a function on $G$ similar to
the ones in the previous sections. Then $I_{\gC_0}(\gph)$ is
finite for any $\gph\in L_2(\gO_0,X)$. Denote the minimum value of
this functional by $M_{\Gamma_{0}}$ and call $\gph$
$\gC_0$--harmonic if $I_{\gC_0}(\gph)=M_{\gC_0}$.

\begin{prop}\label{prop:existcont}
There exists a $\Gamma_{0}$--harmonic map.
\end{prop}

This is actually a special case of Theorem \ref{thm:existence}, where $G$ can be considered as a product $G\times 1$. However, since Proposition \ref{prop:existcont} is much simpler than the general case of Theorem \ref{thm:existence} we shall give an alternative simpler proof. 

\begin{proof}
It is straightforward to check that
\[
I_{\Gamma_{0}}(\gph)=\frac{1}{\mu(\gO_0)}\int_{\Omega_0\times\Omega_0}\sum_{\gamma\in\Gamma_{0}%
}h(\go_{1}^{-1}\gamma
\go_{2})d(\gph(\go_{1}),\gamma\gph(\go_{2}))^{2}.
\]
Let $\gS$ be finite generating set for $\gC_0$. By Lemma
\ref{lem:b>0} $d_\gS (x)\to\infty$ in at least a linear rate.
Since $h\geq\gb>0$ on $\Omega_0(\Sigma\cup\{1\})\Omega_0,$ we have
\[
I_{\Gamma_{0}}(\varphi)\geq\frac{\gb}{\mu
(\gO_0)}\int_{\Omega_0}\left( \int_{\Omega_0}\sum
_{\gamma\in\Sigma\cup\{Id\}}d(\varphi(\go_{1}),\gamma\varphi(\go_{2}))^{2}d\mu
(\go_{2})\right) d\mu(\go_{1}).
\]
We divide the integration in $\go_{2}$ into two parts according to
weather $d(\varphi(\go_{2}),x_0)\geq d(\gph (\go_1),x_0)/2$ or
not. Taking in account only the $\gamma\in\Sigma\cup\{1\}$ which
gives the largest contribution, in view of Lemma \ref{lem:b>0}, we
get for some constant $c>0$
\[
I_{\Gamma_{0}}(\varphi)\geq
c\int_{\Omega}d(\varphi(\go_{1}),x_{0})^{2}=c\left\Vert
\varphi\right\Vert _{2}^{2}.
\]
This means that for any minimizing sequence $\varphi_{n}$ the
$L_{2}$--norm is uniformly bounded. As in the proof of Theorem
\ref{thm:existence}, we can find such a sequence which is cauchy
and hence converges to a harmonic map.
\end{proof}

\begin{prop}
The $\Gamma_{0}$--harmonic map is unique.
\end{prop}

\begin{proof}
Take two $\Gamma_{0}$--harmonic maps $\varphi$ and $\psi.$ BNPC
implies that
\[
d\big(\frac{\varphi(g_{1})+\psi(g_{1})}{2},\frac{\varphi(g_{2})+\psi(g_{2})}{2}%
\big)^2<\frac{
d(\varphi(g_{1}),\varphi(g_{2}))^{2}+d(\psi(g_{1}),\psi
(g_{2}))^{2}}{2}
\]
unless $[\gph (g_1),\psi (g_1)]\parallel [\gph (g_2),\psi (g_2)]$.
Since both $\gph$ and $\psi$ minimize $I_{\gC_0}$ and are
$\gC_0$--equivariant, and since there are no parallel orbits for
the $\gC_0$--action, it follows that $\gph=\psi$.
\end{proof}

\begin{lem}\label{lem:lambda}
Let $\varphi$ be the (unique) $\Gamma_{0}$--harmonic map and
$\lambda\in\Lambda.$ Then
$\varphi_{\lambda}(g):=\lambda^{-1}\varphi(\lambda g)$ is the
$\lambda^{-1}\Gamma _{0}\lambda$--harmonic map.
\end{lem}

\begin{proof}
Indeed
$$
 \varphi_{\lambda}(\lambda^{-1}\gamma\lambda
 g)=\lambda^{-1}\varphi(\gamma\lambda
 g)=\lambda^{-1}\gamma\varphi(\lambda
 g)=\lambda^{-1}\gamma\lambda\varphi_{\lambda }(g),
$$
and hence $\gph_\gl$ is $\lambda^{-1}\Gamma
_{0}\lambda$--equivariant. It is also straightforward to verify
that $I_{\gC_0}(\gph )=I_{\gl^{-1}\gC_0\gl}(\gph_\gl )$. Note that
$\gl^{-1}\gO_0$ is a fundamental domain for $\gl^{-1}\gC_0\gl$.
\end{proof}

\begin{lem}\label{lem:normal}
Let $\Gamma_{1}$ be a finite index normal subgroup of $\Gamma_{0}$
and let $\varphi$ be the $\gC_1$--harmonic map. Then $\varphi$ is
also $\Gamma_{0}$--harmonic.
\end{lem}

\begin{proof}
By normality and Lemma \ref{lem:lambda} we get
$\varphi_{\gamma}(g)=\varphi(g)$ for any $\gamma\in\Gamma_{0}$,
which proves the $\Gamma_{0}$--equivariance. Additionally,
choosing the fundamental domain for $\gC_1$ to be the union of
finitely many translations of $\gO_0$ one can easily verify that
$I_{\gC_1}(\gph')=I_{\gC_1}(\gph')$ for any $\gC_0$--equivariant
map.
\end{proof}

\begin{lem}\label{1234}
Let $\gC_i\leq\gD,i=1,2$ be two subgroups commensurable to $\gC$
with associated harmonic maps $\gph_i$. Then $\gph_1=\gph_2$.
\end{lem}

\begin{proof}
Take $\Gamma_{4}\leq\Gamma_{3}$ of finite index in
$\Gamma_{1}\cap\Gamma_{2}$ such that $\Gamma_{3}$ is normal in
$\Gamma_{1}$ and $\Gamma_{4}$ is normal in $\Gamma_{2}$. Lemma
\ref{lem:normal} implies that $\varphi_{1}=\varphi
_{3}=\varphi_{4}=\varphi_{2}$, where $\gph_i$ is the
$\gC_i$--harmonic map, $1\leq i\leq 4$.
\end{proof}

We conclude:

\begin{prop}
The $\gC$ harmonic map $\gph$ is $\Lambda$--equivariant.
\end{prop}

\begin{proof}
Let $\gl\in\Lambda$. By Lemma \ref{1234} $\gph$ is also
$\gl^{-1}\gC\gl$--harmonic, and hence by Lemma \ref{lem:lambda}
$\gph=\gph_\gl:=\gl^{-1}\gph(\gl\cdot)$. Since $\gl\in\gD$ is
arbitrary, $\gph$ is $\gD$--equivariant.
\end{proof}

Since $\Lambda$ is dense in $G$, and hence acts ergodically on $G$
and since $\gph$ is measurable and $\Lambda$--equivariant, we
conclude:

\begin{cor}
The harmonic map $\gph$ is essentially continuous.
\end{cor}

\begin{proof}
Indeed, for each $g'\in G$ the function $g\mapsto d(\gph (g),\gph (gg'))$ is measurable and $\gD$--invariant, hence constant. The result follows as the action by right translations of $G$ on $L_2(\gO,X)$ is continuous.
\end{proof}

Changing $\gph$ on a set of measure $0$, we can assume that it is
actually continuous. We derive from continuity and
$\Lambda$--equivariance that the set of points $x\in X$ for which
the orbit map $\gl\mapsto\gl\cdot x$ from $\Lambda$ to $X$ is
continuous, is non-empty, indeed it contains $\gph (G)$. Since
this set is also convex and closed, it follows from
$C$--minimality that the orbit map is continuous for every $x\in
X$, and hence that the isometric action extends continuously to
$G$.

\begin{rem}
$(i)$ As in Section \ref{sec:non-uniform}, also Theorem
\ref{thmcommensurator} can be generalized to
the setting of finitely generated ($1<p<\infty$)--integrable
non-uniform weakly cocompact lattices. 
The only part that needs new justification is Proposition \ref{prop:existcont}. However since Proposition 
\ref{prop:existcont} is a special case of Theorem \ref{thm:existence} this generalization can be derived from the discussion in the previous section. 
Moreover, the assumption that $\gC$ is weakly cocompact is actually not required here. One can see this by arguing as in the proof given in Section \ref{sec:existence}, for the special case that $G_2=1$. In order to keep this section as simple as possible we chose to state and give a complete proof of Theorem \ref{thmcommensurator} under the assumption that $\gC$ is uniform and only remark on how this can be generalized to $p$--integrable lattices using arguments that appeared in earlier sections. 
 
$(ii)$ In \cite{Monod} it is not assumed that there are no
parallel orbits, but the conclusion is much weaker: the
$\gC$--action extends but not necessarily the $\Lambda$--action as
in Theorem \ref{thmcommensurator}. The following example shows
that this assumption is required here: 
 Let $G$ be the full
isometry group of $\mathbb{R}$ including the reflection $f$ in
$0$. Let $\Gamma=\langle f\rangle\ltimes\BZ$ and $\Lambda=\langle
f\rangle\ltimes\BZ[\sqrt{2}]$, and let $\Lambda$ act on $\mathbb
R$ where $f$ acts by reflection, $\BZ$ by positive and
$\BZ\cdot\sqrt{2}$ by negative translations.
In this example also the displacement (of the index 2 subgroup $\BZ$ of $\gC$) does not go to infinity (see Footnote \ref{fn:disp-->para}), however, one can still prove that harmonic maps do exists, the lack of uniqueness is what prevents us from extending the $\Lambda$ action to $G$.
\end{rem}

\end{document}